\documentstyle[12pt]{article}

\newlength{\abstractwidth}
\renewcommand{\thefootnote}{\fnsymbol{footnote}}
\renewcommand{\thanks}[1]{\footnote{#1}} 
\newcommand{\starttext}{ \setcounter{footnote}{0}
\renewcommand{\thefootnote}{\arabic{footnote}}}

\newcommand{\be}{\begin{equation}}
\newcommand{\bea}{\begin{eqnarray}}
\newcommand{\eea}{\end{eqnarray}} \newcommand{\ee}{\end{equation}}
\newcommand{\N}{{\cal N}} 
 \def\ba{\begin{eqnarray}}
\def\ea{\end{eqnarray}}

\def\N{{\cal N}}
\def\G{{\cal G}}

\def\C{{\bf C}}

\def\o{\omega}

\def\det{{\rm det}}

\def\log{\,{\rm log}\,}

\def\o{\omega}

\def\a{\alpha}

\def\b{\beta}
\def\g{\gamma}
\def\d{\delta}

\def\o{\omega}

\def\O{\Omega}

\def\na{\nabla}

\def\N{\bf N}

\def\R{{\bf R}}
\def\C{{\bf C}}

\def\i{\infty}

\def\p{\partial}

\def\G{\Gamma}

\def\na{{\nabla}}

\def\G{\Gamma}

\def\[{{\bf [}}
\def\]{{\bf ]}}

{
	\par\vspace{\baselineskip}%
	\textbf{Proof}\begin{itshape}%
		\par\vspace{\baselineskip}%
	}%
	{\end{itshape}}

\begin{document}
\starttext \baselineskip=18pt \setcounter{footnote}{0}
\newtheorem{theorem}{Theorem}
\newtheorem{lemma}{Lemma}
\newtheorem{corollary}{Corollary}
\newtheorem{definition}{Definition}
\newtheorem{conjecture}{Conjecture}
\newtheorem{proposition}{Proposition}
\newtheorem{example}{Example}

\begin{center}
{\Large \bf Convergence of Hermitian manifolds and the Type IIB flow
}

\medskip
\centerline{Nikita Klemyatin\footnote{Supported in part by the National Science Foundation under grant DMS-1855947.}}

\medskip

\begin{abstract}

{\footnotesize The Type IIB flow is a flow of conformally balanced complex manifolds introduced by Phong, Picard, and Zhang, about whose singularities little is as yet known. We formulate convergence criteria for the Gromov-Cheeger-Hamilton convergence of sequences of Hermitian manifolds, and apply them to precompactness theorems and the existence of singularity models for the Type IIB flow, in analogy with Hamilton's classic compactness theorems and classification of singularities for the Ricci flow.}

\end{abstract}

\end{center}

\baselineskip=15pt
\setcounter{equation}{0}
\setcounter{footnote}{0}


\section{Introduction}
\setcounter{equation}{0}

The earliest, and still most influential, geometric structure which emerged as a candidate for vacuum configuration of unified string theories was that of a Calabi-Yau manifold, equipped with a holomorphic vector bundle admitting a Hermitian-Einstein metric. This was proposed by Candelas, Horowitz, Strominger, and Witten \cite{CHSW} for the heterotic string, and it could rely on the fundamental existence theorems for such structures proved earlier by Yau \cite{Y} and Donaldson-Uhlenbeck-Yau \cite{D,UY}. But other geometric structures have emerged since, motivated by the other string theories, for which analogous existence theorems are not yet available. Such theorems can be expected to require the full theory of non-linear partial differential equations and, as advocated in \cite{PPZ17, P}, notably the theory of geometric flows. In particular, flow approaches have been proposed in \cite{PPZ17} for the Hull-Strominger system in heterotic string theory, in \cite{P,FPPZa} for the Type IIB string, in \cite{FPPZb} for the Type IIA string, and in \cite{FGP} for $11$-dimensional supergravity.

\smallskip
These new flows share main features with the Ricci flow, with typically some additional structure, but they remain more complicated than, say, the K\"ahler-Ricci flow. A few known cases of convergence are in \cite{PPZAnn, PPZProc, PPZMathAnn}, but otherwise very little is known as yet about their possible singularities. It is a natural first step to establish the existence of singularity models, in analogy with the case of the Ricci flow. In this paper, we carry this out for the Type IIB flow.

\smallskip
We state now our main results, The definitions of all relevant notions are given in \S 2 below. The first two results are analogues for non-K\"ahler Calabi-Yau manifolds and for the Type IIB flow of Hamilton's compactness theorems for the Ricci flow on Riemannian manifolds \cite{Hamilton1}:
\begin{theorem}\label{CompactnessTheorem}
\label{compactness}
	Let $\{(M_j, g_j, J_j,\Omega_j, x_j)\}_{j=1}^\infty$ be the sequence of compact marked Hermitian Calabi-Yau manifolds. Assume that we have the following properties:
	
	\smallskip

{\rm (1)} The sequence has bounded Hermitian geometry, in the sense of Definition \ref{BoundHermGeom} below;

{\rm (2)} There are constants $C_1 \geq C_0 > -\infty$, such that $C_1\geq\log\|\Omega_j\|^2_{g_j}\geq C_0$.
		
		
		{\rm (3)} All manifolds are conformally balanced: $d(\|\Omega_j\|_{g_j}\omega^{n-1}_j)=0$.
		
		\smallskip
	
	Then, after passing possibly to a subsequence, there is a complete marked Hermitian Calabi-Yau manifold $(M_\i, g_\i,J_\i, \Omega_\i, x_\i)$, an exhaustion $\{U_j\}_{j \in \N}$ of $M_\i$ and a sequence of diffeomorphisms $\Phi_j: U_j \rightarrow V_j \subset M_j$, $\Phi_j(x_\i) = x_j$, such that the triple $(\Phi_j^*g_j, \Phi_j^*J_j, \Phi_j\Omega_j)$ converge to the limit triple $(g_\infty, J_\infty, \Omega_\infty)$. The limit complex structure is integrable and $\Omega_\infty$ is a holomorphic $(n,0)$-form with respect to it. Moreover, the limit metric is conformally balanced. 
\end{theorem}

\smallskip

\begin{theorem}
\label{compactnessTypeIIB}
	Let $\{(M_j, g_j(t), J_j,  \O_j, x_j)\}_{j=1}^{\infty}$ be a sequence of solutions of the Type IIB flow on the interval $[\tau_0; \tau_1]$, such that:
	
	\smallskip
	
		{\rm (1)} There is a constant $C_0$, such that $$|Rm_j|_{g_j} + |\na T|_{g_j} + |T|^2_{g_j}\leq C_0$$
		
		
		{\rm (2)} There is a constant $\iota_0>0$, such that $\mathrm{inj}(M_j) \geq \iota_0$;
		
		{\rm (3)}There is a constant $C_1$, independent of $j$, such that $C_1 \leq \log||\Omega_j||^2_{g_j}(\tau_0)$;
		
		\smallskip

	Then, after passing possibly to a subsequence, there is a complete conformally balanced manifold $(M_\infty, g_\infty(t), J_\infty, \Omega_\infty, x_\infty)$, an exhaustion $\{U_j \}_{j \in \N}$ of $M_\i$, and a sequence of diffeomorphisms $\Phi_j: U_j \rightarrow V_j \subset M_j$, such that $(M_j, g_j(t), J_j, \O_j, x_j)$ converge to a complete solution $(M_\infty, g_\infty(t), J_\infty, \Omega_\infty, x_\infty)$ of the Type IIB flow. 
\end{theorem}

Our third result is the existence of singularity models for the Type IIB flow, in analogy with Hamilton's theorem for the Ricci flow \cite{H95}. A singularity model is a solution of the Type IIB flow on a (possibly non-compact) conformally balanced complex manifold, which arises as a blow-up limit of a Type IIB flow, as the time parameter $t$ approaches its maximum time of existence. 

\begin{theorem}
\label{models}
Let $\o(t)$ be a solution of the Type IIB flow on a compact conformally balanced manifold $M$, and let $[0,T)$ be its maximum time of existence. Assume that the solution admits an injectivity radius estimate, in the sense of Definition \ref{injestimate} below. Then for suitable times $t_j\to T$ and suitable rescalings $g_j(t)=C_jg(t_j+C_j^{-1}t)$, the sequence $g_j(t)$ admits a convergent subsequence to a solution of the Type IIB flow with a corresponding function $f(t)$, as defined in equation (\ref{f}) below, satisfying the following properties:

\smallskip
	
		{\rm (1)}  For the Type I the solution exists on the interval $(-\infty, c)$ and $f(t) \leq {C \over c - t}$;
		
		{\rm (2)} For the Type IIa the solution exists for all time and $f(t) \leq C$;
		
	{\rm (3)} For the Type IIb the solution exists for all time and $f(t) \leq C$ (the same as the previous case);
	
		{\rm (4)} For the Type III singular solution the singular model solution exists on the interval $(-A, \infty)$ and $f(t) \leq {A \over A + t}$.
	
\end{theorem}

\medskip
Finally, we point out a new way of expressing the Type IIB flow, which reveals its similarities with Bryant's $G_2$ Laplacian flow
\cite{B,BX} as well as the Laplacian formulation of the Type IIA flow \cite{FPPZb}:

\begin{theorem}
\label{laplacianth}
Let $\langle \alpha, \beta\rangle_\omega$ the standard pointwise Hermitian product of forms and consider the following $L^2$ scalar product on forms: 
\bea
(\alpha, \beta) = \int_M \langle \alpha, \beta\rangle_\omega i^{n^2} \Omega \wedge \overline{\Omega} = \int_M \alpha \wedge \|\Omega\|^2_\omega \overline{*\beta}. 
\eea
Let $\partial_{\o,\Omega}^\dagger$ the operator conjugated to $\partial$ with respect to $(\alpha, \beta)$. Then the conformally balanced condition for $\omega$ is equivalent to the condition $\p_{\o,\Omega}^\dagger \o=0$, and the Type IIB flow, as defined originally in (\ref{TypeIIB-omega}) below, can be rewritten as
\bea
\label{laplacian}
\partial_t \omega = - \partial_{\o,\Omega}^\dagger\partial\omega.
\eea
\end{theorem}

{\bf Acknowledgments:} I am indebted to Professor Duong H. Phong for suggesting this problem and his constant attention and help while the paper was in preparation. I am also grateful to Misha Verbitsky for the fruitful discussions about the Gromov-Hausdorff convergence.

I was partially supported by the HSE University Basic Research Program, Russian Academic Excellence Project '5-100', by the M\"obius Contest Foundation for Young Scientists, and by the "Young Russian Mathematics" award.

\section{Preliminaries}
\setcounter{equation}{0}

We begin by describing our notations and conventions. 
Let $(M, g, J)$ be a Hermitian manifold, and $\omega$ be the K\"ahler form:, $\omega(X,Y) = g(X,JY)$. In local holomorphic coordinates, we write
$$\omega = ig_{j\bar{k}}dz^j\wedge d\bar{z}^k $$
Let $V = V^j \partial_j$ be a section of $T^{1,0}M$. The Chern connection $\na$ is given by
\bea
&&
\nabla_j V^k = \partial_jV^k + \Gamma^k_{js} V^s, 
\quad
\nabla_{\bar{j}} V^k = \partial_{\bar{j}}V^k,
\nonumber\\
&&
\Gamma^k_{js} = g^{k\bar{p}}\partial_jg_{s\bar{p}}.
\eea
The curvature tensor $Rm=R_{\bar{k} j}{}^p{}_q$ is then given by
\bea
R_{\bar{k} j}{}^p{}_q = -\partial_{\bar{k}}\Gamma^p_{jq} 
\eea
and the first and second Chern-Ricci curvatures are defined by
\bea
&&
R_{\bar{k}j} = R_{\bar{k} j}{}^p{}_p = - \partial_{\bar{k}}\partial_j\log\det(g) \nonumber\\
&&
\tilde{R}_{\bar{p} q} = g_{\bar{p} \ell} g^{j \bar{k}} R_{\bar{k} j}{}^\ell{}_q.
\eea
The torsion is defined by 
\bea
T=i\partial \omega, ~ \bar{T} = -i\bar{\partial}\omega.
\eea
If we express $T$ in components as $T = {1 \over 2}T_{\bar l k j}dz^j \wedge dz^k \wedge d\bar{z}^l$,
then $T_{\bar l k j} = \partial_kg_{\bar l j} - \partial_jg_{\bar l j}$.
We also define the contracted torsion $\tau = \Lambda i\partial\omega$ which is a $(1,0)$-form.
In coordinates, $\tau=\tau_jdz^j$, with $\tau_j = g^{k\bar l}T_{\bar l k j}$. 

\subsection{Connections on Hermitian manifolds}

For our purposes, the Chern connection may not be the most convenient Hermitian connection to use. The reason is that we would like to make use of convergence theorems of Riemannian manifolds. These require lower bounds on the injectivity radius with respect to the Levi-Civita connection, which is different from the injectivity radius defined with respect to the Chern connection.

\smallskip
Recall that, given a complex structure $J$ and a Hermitian metric on the manifold, the Chen connection is only one connection in a whole line of Hermitian connections $\nabla^t$, known as the Gauduchon line, and given explicitly by \cite{Gau, AOUV}, 
\bea
\label{GauduchonConnections}
g(\nabla^s_X Y,Z) = g(\nabla^g_X Y,Z) + \frac{1 -t}{4}d^c\omega(X,Y,Z) + \frac{1 + t}{4}d\omega(JX,Y,Z).
\eea
Here $\nabla^g$ denotes the Levi-Civita connection. We shall however be interested in the Bismut connection $\na^+$, which corresponds to $t=-1$. Explicitly, it is given by
\bea
&&
\nabla^+_j V^k = \partial_jV^k + \Gamma^k_{mj}V^m = \nabla_j V^k + T^k_{mj}V^m
\nonumber\\
&&
\nabla^+_{\bar{j}}V^k = \partial_{\bar{j}}V^k + g^{k\bar{m}}\bar{T}_{p \bar{j}\bar{m}}V^p.
\eea
Using the Chern and the Bismut connections as reference points, one can write the Gauduchon line as follows:
\bea
\na^t = (1 - \kappa)\na + \kappa\na^+,
\qquad
\kappa={1-t\over 2}.
\eea

We are interested in the Bismut connection because of  its following property:

\begin{lemma}\label{BismutGeodesics}
The notions of geodesic with respect to the Levi-Civita connection and with respect to the Bismut connection coincide. In particular, both connections define the same notion of injectivity radius.
\end{lemma}

\noindent
{\it Proof.} By definition, the geodesics for the Bismut connection are those curves $x(t)$, which satisfy the equation $\nabla^+_{\dot{x}}\dot{x} = 0$. For any vector field $Z$ we have 
\bea
g(\nabla^+_{\dot{x}} \dot{x},Z) 
= 
g(\nabla^{g}_{\dot{x}} \dot{x},Z) + \frac{1}{2}d^c\omega(\dot{x},\dot{x},Z) .
\eea
But $d^c\omega(\dot{x},\dot{x},Z) = (-Jd\omega)(\dot{x},\dot{x},Z) = 0.$
Hence $\nabla^+_{\dot{x}} \dot{x} = 0$ iff $\nabla^{g}_{\dot{x}} \dot{x}=0.$
The lemma is proved. Q.E.D.

\smallskip
In view of this lemma, we can now introduce a notion of bounded Hermitian geometry formulated entirely in terms of the Bismut connection:

\begin{definition}\label{BoundHermGeom}
	Let $\{(M_j,g_j,J_j)\}_{j \in \N}$ be a sequence of complete Hermitian manifolds. Denote by $Rm^+_j$, and $T^+_j$ the curvature and torsion of the Bismut connection on $(M_j,g_j,J_j)$. We say that this sequence has bounded Hermitian geometry, if the following holds:
	
	\smallskip
		
       {\rm (1)} There is a constant $\iota_0>0$, such that the injectivity radii of these manifolds are bounded below by $\iota_0$;
		
{\rm (2)} There are constants $C_m>0$, independent on $j$, such that 
		\bea
		|(\na^+)^mRm^+|_{g_j} \leq C_m,\qquad
		|(\na^+)^mT^+|_{g_j} \leq C_m
		\eea 
		for any $m \in \N$.
	
\end{definition}

\subsection{Basic facts about Hermitian Calabi-Yau manifolds}

In this section, we collect some facts about the non-K\"ahler Calabi-Yau manifolds. Recall that a 
Hermitian Calabi-Yau manifold is a 4-tuple $(M,g,J,\Omega)$, where $(M,g,J)$ is a Hermitian manifold, and $\Omega$ is a nonvanishing $(n,0)$-form.

\medskip
A Hermitian Calabi-Yau manifold does not always possess a K\"ahler metric. There are many examples of such manifolds: Calabi-Yau manifolds obtained from the conifold transition (see \cite{_CPY_}), hyperk\"ahler fibrations over Riemannian surfaces (\cite{_Fei16_}, \cite{_Fei15_}), factors of Lie groups, etc. A less stringent notion than K\"ahler is that of balanced metric, introduced by Michelson 
\cite{MI82}. These are metrics satisfying $d\omega^{n-1}=0$. Closely related is the notion of conformally balanced metric, defined by the condition 
 $d(\|\Omega\|_\omega\omega^{n-1}) = 0$. The conformally balanced condition was recognized by Li and Yau \cite{LY} as one of the equations in the Hull-Strominger system in string theory. In terms of the metric $\eta = \|\Omega\|\omega$, it can also be expressed as $d(\|\Omega\|^2 \eta^{n-1})= 0$.

\subsection{Pointed Gromov-Cheeger-Hamilton convergence}

We recall the notion of convergence due to Hamilton \cite{Hamilton1}, which is an adaptation of a notion of $C^\infty$ convergence of manifolds due to Gromov and Cheeger:

\begin{definition}
	Let $\{(M_j, g_j, O_j)\}$ be a sequence of complete Riemannian manifolds, where $O_j \in M_j$ is a point.
	We say that this sequence converges to the complete pointed manifold $(M_\infty, g_\infty, O_\infty, q_\infty)$ in $C^k$-sense, if there is a sequence of open sets $\{U_j\}$ in $M_\infty$ with the following properties:

{\rm (1)} The sequence $\{U_j\}$ is an exhaustion of $M_\infty$;

{\rm (2)} There are diffeomorphisms $\Phi_j: U_j \rightarrow V_j \subset M_j$ such that $\Phi_j(O_\infty) = O_j$ and the sequence of metrics $\{\Phi_j^* g_j\}$ converges to $g_\infty$ uniformly in $C^k$-topology on compact subsets in $U_j$. 
	
\end{definition}

The following theorem of Hamilton \cite[Theorem 1.2]{Hamilton1} gives an important set of conditions implying the existence of a convergent subsequence: 
\hfill

\begin{theorem}\label{HamiltonCompactness}
	Let $\{M_j, g_j,O_j\}$ be a sequence of complete pointed Riemannian manifolds and $\nabla_j$ be the Levi-Civita connection on $(M_j,g_j)$. Assume that the following holds uniformly in $j$:
	
	{\rm (1)} There are constants $C_k>0$, such that $|\nabla^k_jRm_j|\leq C_k$;
	
	{\rm (2)} There is a positive uniform lower bound $\iota_0$ for the injectivity radii of $M_j$ at $O_j$, $inj(M_j,O_j)\geq \iota_0$.
	
	Then there is a subsequence of the initial sequence, which converges to a complete pointed Riemannian manifold $(M_\infty, g_\infty, O_\infty)$.
\end{theorem}

\subsection{Curvature estimates}

In this section, we fix a smooth compact manifold $M$ with metric $g_{ij}$. Let $E$ be a smooth vector bundle over $M$ equipped with a metric $H_{\alpha\beta}$. Let $A_1,A_2$ be two connections on $E$, and $F_1,F_2$ be their curvatures. We are interested in obtaining bounds for the curvature of one of them from bounds for the curvature of the other.

\smallskip
Let $C=A_1 - A_2$, which is  a $End(E)$-valued 1-form. We have
\bea
F_1 &= &
dA_1 + F_1 \wedge A_1 
=  
dA_2 + A_2\wedge A_2 + dC + [A_2,C] + C\wedge C 
\nonumber\\
&= &
F_2 + d_{A_2}C + C \wedge C.
\eea
Here $d_{A_2}C$ is the covariant derivative with respect to $A_2$. Let $x^1, \dots, x^n$ be a system of local coordinates on $M$.
Denote by $\nabla^j$ the covariant derivative with respect to $A_j$. Then for any section $\psi$ of $E$ we have
\bea
\nabla^j_a \psi 
= 
\partial_a \psi + (A_j)_a\psi.
\eea 
For our endomorphism-valued 1-form $C = C_adx^a$ we obtain the following formula:
\bea
\nabla^j_aC_b 
= 
\partial_jC_b + [(A_j)_a,C_b].
\eea
In this notation, $d_{A_j}C = (\nabla^j_aC_b - \nabla^j_b C_a)dx^a \wedge dx^b$. Let $F_j = (F_j)_{ab}dx^a\wedge dx^b$ for $j=1,2$. Combining all this, we obtain the following formula:
\bea
(F_1)_{ab} 
= 
(F_2)_{ab} + (\nabla^2_aC_b - \nabla^2_b C_a) + [C_a, C_b].
\eea

Assume now that we have bounds for all covariant derivatives of $F_2$ and $C$ with respect to $\nabla^2$. More presicely, for any natural $p$, there are constants $C_p>0$, such that $|(\nabla^2)^pF_2| < C_p$ and $|(\nabla^2)^pC| < C_p$. We want to establish the estimates of $|(\nabla^1)^pF_1|$ in terms of $|(\nabla^2)^qC|$ and $|(\nabla^2)^pF_2|$.

\smallskip

In order to do this, we recall some basic facts about the norms of tensors. First of all, for any two $E$-valued tensors $T_1, T_2$ on $M$, we have a pointwise inequality:
\bea
|T_1 * T_2| \leq C_1|T_1||T_2| 
\eea
where $\star$ denotes any bilinear operation involving the metrics $g_{ij}$ and $H_{\alpha\beta}$, and the constant $K_1=K_1(dim(M), rk(E))$ is a universal constant depends only on the dimension $n$ of $M$, the rank of $E$, and the type of $T_1$ and $T_2$, but independent of $g_{ij}$ and $H_{\alpha \beta}$. This can be proven by Linear Algebra.

Next, for any tensor $T$ we have the following 
\bea
\label{VeryImportantIneq}
|(\nabla^2 + C)T| 
\leq 
|\nabla^2 T| + |CT|
\leq 
|\nabla^2 T| + C_1|C||T|. 
\eea
Iterating this process, we obtain the following inequality:
\bea
\label{HighDeriv}
|(\nabla^2 + C)^pT|\leq |(\nabla^2)^p T|+\Sigma_{k=0}^{p-1}|C*\dots*(\nabla^2)^k(C*\dots*C*T)| &\leq&
\nonumber\\  
|(\nabla^2)^p T| + K_p\Sigma_{k_1 + \dots + k_p = p-1}|(\nabla^2)^{k_1}C|\dots|(\nabla^2)^{k_{p-1}}C||(\nabla^2)^{k_p}T|,
\eea
where $K_p$ are universal, and independent of the metrics $g_{ij}$ and $H_{\alpha \beta}$
Now we can state the main lemma of this section.

\begin{lemma}\label{CurvatureEstimates}
	Let $M$ be a smooth manifold, E is a vector bundle over $M$, $A_1, A_2$ are arbitrary connections and $F_1, F_2$ are they curvatures. Let $C=A_1 - A_2$ (as above). Suppose that there are constants $C_m>0$, such that $|(\na^2)^mF_2| < C_m$ and $|(\na^2)^mC| < C_m$ for any natural $m$. Then there are constants $K_m$ depending on $C_p$ and independent on the metrics $g_{ij}$ and $H_{\alpha\beta}$ such that $|(\na^1)^mF_1|<K_m$
\end{lemma}
The proof is just applying the inequality above to $T=F_2$, $T=d_{A_2}C$ and $T=[C,C]$. Note that there is a constant $C_{1,1}$, such that $|[C,C]|< C_{1,1}|C|^2$. This lemma has very important consequences:

\begin{lemma}\label{RiemannianEstiamtes}
	Let $(M,g,J)$ be a Hermitian manifold. Let  $\na$ and $\na^g$ be the Chern and the Levi-Civita connections on $M$ with curvatures $Rm, Rm^g$. Suppose there are constants $C_m$, $m \in \N$, such that $|\na^m Rm|_g \leq C_m$ and $|\na^m T|_g \leq C_m$. Then there are constants $K_m$, such that $|(\na^g)^m Rm^g|_g\leq K_m$.
\end{lemma}

\noindent{\em Proof.}
Let $C(X,Y) = \na_X Y - \na^g_X Y$ for any vector fields $X,Y$ on $M$. This is a bilinear form with values in the tangent bundle. By the formula \ref{GauduchonConnections} , we have the following:
\bea
g(C(X,Y),Z)
&=&
g(\na_X Y,Z) - g(\na^g_X Y,Z) 
= 
d\o(JX,Y,Z)
\nonumber\\
&=&
(T^{\R}*J*J)(X,Y,Z)
\eea
where $T^\R = T + \bar T = i(\p - \bar \p)\o$. Hence, there is a constant $P$, such that $|C|_g \leq P|T|_g|J|_g$. Similarly, $|\na^m C|_g = |J * \na^mT|_g \leq P_m |\na^mT|_g|J|_g$ for some universal constant $P_m$. Note that $|J|_g = \sqrt{n}$, where $n = \dim_\C M$.
Hence, we obtain bounds for $|\na^m C|_g$. Now we apply Lemma \ref{CurvatureEstimates} with $\na^2 = \na$, $\na^1 = \na^g$. Q.E.D. 

\begin{lemma}\label{BoundsOnComplexStructure}
	Let $(M,g(t),J)$ be a family Hermitian manifold, $t \in \R$. Assume that $|\na^mT|_{g(t)} \leq C_m$, where $C_m$ are independent on $g(t)$. Then there are constants $K_m$, such that $|(\na^g)^mJ| \leq K_m$, and these constants are independent on $g(t)$.
\end{lemma}

\noindent{\em Proof.}
For  notational simplicity, we omit the dependence on $t$ of $g(t)$. By definition of the Chern connection, $\na J = 0$. Hence,
\bea
\na^g J
=
\na^g J - \na J
=
T^\R * J * J * J
\eea
and
\bea
|\na^g J|_g \leq P_1 |T|_g |J|^3_g.
\eea
The constant $P_1$ depends only on the dimension $n=\dim_\C M$. The norm $|J|_g$ is equal to $\sqrt{n}$, hence 
\bea
|\na^g J|_g \leq P_1 C_1 n^{3 \over 2}.
\eea
In particular, this bound does not depend on the metric. For the case of higher derivatives we can use the formula (\ref{HighDeriv}). We have
\bea
|(\na^g)^mJ|_g
\leq
P_m\sum_{k_1 + \dots + k_{m-1} = m-1}|\na^{k_1}T|_g \dots |\na^{k_{m-1}}T|_g|J|^{q(m)}_g,
\eea
and $P_m, q(m)$ does not depend on metric. Since the all derivatives of $T$ are bounded, we conclude, what there are constants $K_m = P_m n^{q(m) \over 2}\sum_{k_1 + \dots + k_{m-1} = m-1}C_{k_1} \dots C_{k_{m-1}} $
 Q.E.D.
\begin{lemma}
	Let $\{(M_j, g_j(t), J_j)\}_{j \in \N}$ be the sequence of Hermitian manifolds with bounded geometry in the sense of the Definition \ref{BoundHermGeom}. Denote by $RM_{g_j}$ the curvature of the Levi-Civita connection on $(M_j,g_j)$.
	Then there are constants $K_m>0$, independent of $j$, such that $|(\na^{g_j})^mRm_{g_j}|_{g_j} \leq K_m$.
\end{lemma}

\noindent{\em Proof.}
In general, the Bismut connection and the Levi-Civita connection are related by the following formula:
\bea
\na^{+} = \na^{g} + {1 \over 2}g^{-1}H
\eea
and $H = d^c\o$. On the other hand, $H = g(T^+(X,Y),Z)$. Applying to the case at hand, and omitting the dependence on $t$ from the notation, we can write $\na^{+,j} - \na^{g} = {1 \over 2}g^{-1}_jH_j$.
The desired statement follows now from Lemma \ref{CurvatureEstimates} with $\na^1 = \na^{g_j}$ and $\na^2 = \na^{+,j}$. Q.E.D.

\begin{lemma}\label{GauduchonEstimates}
	Let $\{(M_j, g_j(t), J_j)\}_{j \in \N}$ be a sequence of Hermitian manifolds, and $\na^j, Rm_j, T_j$ be the Chern connection, its curvature and torsion respectively. Denote by $\na^{s,j}, Rm_{s,j}, T_{s,j}$ be the Gauduchon connection, its curvature and torsion with parameter $s$ on the Gauduchon line. Assume that there are constants $C_m>0$, independent on $j$, such that $|(\na^j)^mRm_j|_{g_j} \leq C_m$ and $|(\na^j)^mT_j|_{g_j}\leq C_m$. Then for any $s\in \R$ there are constants $K_{s,m} >0$, independent of $j$, such that 
	\bea
	|(\na^{s,j})^m Rm_{t,j}|_{g_j} \leq K_{s,m},\qquad
	|(\na^{s,j})^m T_{t,j}|_{g_j} \leq K_{s,m}.
	\eea
\end{lemma}

\noindent{\em Proof.}
We know that $\na^s = \na^g + {1 + s\over 4}g^{-1}H + {1 - s \over 4}g^{-1}S$, where $H(X,Y,Z) = d^c\o (X,Y,Z)$ and $S(X,Y,Z) = d\o(JX,Y,Z)$. Next, $\na^g = \na - {1 \over 2}g^{-1}S$. Hence
\bea
\na^s = \na + {1 + s\over 4}g^{-1}H - {1 + s \over 4}g^{-1}S.
\eea 
This implies that the difference $\na^t - \na$ has the form of $T^\R * J * J$ and $T^\R$, where $T^R=d^c\o$. The bounds on $T_j$ implies the bounds on $T^\R_j$ with all derivatives. Thus, the difference $\na^t - \na$ has bounded covariant derivatives w.r.t. $\na^j$.  Lemma \ref{CurvatureEstimates} implies bounds for $|(\na^{t,j})^m Rm_{t,j}|_{g_j}$. The bounds for $|(\na^{t,j})^m T_{t,j}|_{g_j}$ follow from the inequality (\ref{VeryImportantIneq}) with $\na^2 = \na^j$ and $\na^1 = \na^{t,j}$. Q.E.D.

\medskip

Finally, we can prove an analogue of Lemma \ref{BoundsOnComplexStructure} for any Gauduchon connection.

\begin{lemma}
	Let $(M,g(t),J)$ be a family of Hermitian manifolds, $s\in \R$, $\na^s$ as above. Assume that $|(\na^s)^mT_s|_{g(t)} \leq C_m$, where $C_m$ are independent of $g(t)$. Then there are constants $K_m$, such that $|(\na^g)^mJ| \leq K_m$, and these constants are independent of $g(t)$.
\end{lemma}

\noindent{\em Proof.}
The proof is similar to the proof of Lemma \ref{BoundsOnComplexStructure}. We have the formula $\na^s = \na^g + {1 + s\over 4}g^{-1}H + {1 - s\over 4}g^{-1}S$, and $\na^sJ = 0$. The bound on torsion implies the bound on $d \o$ by the following formula
\bea
(d \o)_{ k jl} 
= 
\na^s_j \o_{ k l} + \na^s_k \o_{l j} + \na^s_l\o_{j k} + 3(T_s*\o)_{kjl}
=
3(T_s*\o)_{kjl},
\eea
where our indices are {\em real} indices.
The same is true for $d^c \o$. Hence we obtain $d^c\o$, and both $H = d^c\o$ and $S = d\o *J$. This implies a bound on the difference $C_s=\na^g - \na^s$ in terms of norms of $H$ and $S$, which are bounded by the norm of the torsion multiplied by a universal constant. Now we apply the inequality (\ref{HighDeriv}) for $(\na^g)^mJ$ with $\na^2 = \na^t$ and $C = C_s$.

\subsection{The Type IIB flow}

Let $M$ be a complex manifold of complex dimension $n$ with trivial canonical bundle $K_M = \det(T^{(1,0)}M)^*$, and let $\Omega$ be a nowhere vanishing holomorphic $(n,0)$-form. 
The Type IIB flow is the flow of $(n-1,n-1)$-forms defined as follows \cite{P}
\bea
\partial_t(\|\Omega\|_\eta \eta^{n-1}) = i\partial \overline{\partial}\eta^{n-2}.
\eea
with conformally balanced intiial data $\o_0$.
It can also be viewed as a special case with vanishing slope parameter $\alpha'$ of the Anomaly flow introduced in \cite{PPZ17}. If we define the metric $\omega$ by $\omega = \|\Omega\|_\eta \eta = ig_{j\bar{k}}dz^j \wedge d\bar{z}^k$, then the above flow for $\eta$ is equivalent to the following flow for $\o$,
\bea
\label{TypeIIB-omega}
\partial_t(\|\Omega\|^2_\omega \omega^{n-1}) = i\partial \overline{\partial}(\|\Omega\|^2_\omega \omega^{n-2}) 
\eea
A motivating property for the Type IIB flow is the fact that it preserves the conformally balanced property of the initial data. For estimates, it is useful to express it as follows as a flow of $(1,1)$-forms 
(see \cite{_FP19_} for the detailed derivation)
\bea
\partial_t g_{\bar{k}j}  = - \tilde{R}_{\bar{k}j} - {1 \over 2}T_{\bar{k}pq}\bar{T}_{j}^{~pq}
\eea
where $\tilde{R}_{j\bar{k}} = g_{s\bar{k}}g^{p\bar{q}}R^s_{~j\bar{q}p}$ is the second Chern-Ricci curvature of the metric. 
The Type IIB flow reduces to the K\"ahler-Ricci flow on K\"ahler manifolds and one can think about it as a particular generalization of the K\"ahler-Ricci flow to non-K\"ahler manifolds
\footnote{Other generalizations have been considered by Streets and Tian \cite{ST} and Ustinovsky \cite{U}.}.
In particular, the Type IIB flow shares many properties with the K\"ahler-Ricci flow. 
As described in Theorem \ref{laplacianth}, we can provide yet another formulation of the Type IIB flow in the spirit of Laplacian flows:

\medskip
\noindent{\it Proof of Theorem \ref{laplacianth}.}
Recall that for any (2,1)-form $\psi$ the usual operator $\partial^\dagger$ acts as follows (see \cite[Appendix D]{_FP19_}):
\bea
(\p^\dagger\psi)_{\bar\alpha\beta}
=
-g^{\gamma\bar j}\na_{\bar j}\psi_{\bar\alpha\beta\gamma}
+
g^{\gamma\bar j}\bar\tau_{\bar j}\psi_{\bar\alpha\beta\gamma}
-
{1\over 2}
\bar T_{\beta\bar j\bar m}\psi_{\bar\alpha\gamma\delta}g^{\gamma\bar j}g^{\delta\bar m}
\eea
The second summand appears from the integration by parts formula for the divergence of the Chern connection:
\bea
\int_M \nabla_jV^j \o^n
=
\int_M \tau_j V^j\o^n
\eea
for any vector field $V^j$. On the other hand, by the conformally balanced condition, we have $\nabla_\|\Omega\|^2_\o = \|\Omega\|^2_\o\tau_j$ (see \cite[Lemma 4]{PPZ16}). This implies 
\bea
\int_M \nabla_jV^j i^{n^2} \Omega \wedge \overline{\Omega}
=
\int_M (\nabla_jV^j)||\Omega||^2_\o  {\o^n \over n!}
=
0.
\eea
Hence, by the same argument as in \cite[Appendix D]{_FP19_}, we can show that for any (2,1)-form $\psi$ we have
\bea
(\p^\dagger_{\o, \O}\psi)_{\bar\alpha\beta}
=
-g^{\gamma\bar j}\na_{\bar j}\psi_{\bar\alpha\beta\gamma}
-
{1\over 2}
\bar T_{\beta\bar j\bar m}\psi_{\bar\alpha\gamma\delta}g^{\gamma\bar j}g^{\delta\bar m}.
\eea
In fact, $\p^\dagger_{\o, \O} = \p^\dagger - \iota_{\tau}$, where $\iota_{\tau}$ is the contraction with the vector field $g^{j\bar k}\bar{\tau}_k$.
The equivalence between the conformally balanced condition and the condition $\p_{\o,\Omega}^\dagger\o=0$ follows.
Next, let $\psi = T = i\p \o$. Then 
\bea
(\p^\dagger_{\o, \O}T)_{\bar\alpha\beta}
=
-\na^\gamma T_{\bar\alpha\beta\gamma} 
-
{1 \over 2}
\bar T_{\beta\bar j\bar m}T_{\bar\alpha\gamma\delta}g^{\gamma\bar j}g^{\delta\bar m}
=
- \tilde{R}_{\bar \alpha \beta} 
- {1 \over 2}T_{\bar\alpha}^{\bar j \bar m}\bar{T}_{\beta\bar j\bar m}
\eea
where in the second equality we used the fact $\tilde{R}_{\bar \alpha \b} = -\na^\gamma T_{\bar \a \b \gamma}$ (see \cite[Lemma 5,iii]{PPZ16}).
Combining all together, and using the fact that $\omega = i g_{\bar k j}dz^j\wedge d\bar{z}^k$, we obtain that
\bea
\p_t(ig_{\bar k j}) 
= 
i(\p^\dagger_{\o,\O}i\p \o)_{\bar k j}
=
-(\p^\dagger_{\o,\O}\p \o)_{\bar k j}.
\eea
Q.E.D.


\medskip

We observe that the evolution equation for the torsion was derived in
\cite[Equation 4.6]{FPPZa}. Using our Laplacian formulation, we can readily rederive this result:
\bea
\p_t T 
= 
\p_t i \p \o
=
- i \p \p^\dagger_{\o, \O} \p \o 
=
-\p \p^\dagger_{\o, \O} T.
\eea

\subsection{Shi-type estimates for connections on the Gauduchon line}

Shi-type estimates for the Type IIB flow were established in  \cite{FPPZa,PPZ16}:

\begin{theorem}
\label{ShiTypePhongEtAl}
Let $(M,\omega)$ be a Calabi-Yau manifold and $\omega$ evolve by the Type IIB flow. Suppose there is a constant $A$, such that $|Rm| + |\nabla T| + |T|^2 \leq A$ for $t \in [0;{1 \over A}]$. Then for any $k > 0$ there is $C_k >0$ (depending on the lower bound for $||\Omega||$), such that $|\nabla^kRm|\leq \frac{C_kA}{t^{k/2}}$ and $|\nabla^{k+1}T|\leq \frac{C_kA}{t^{k/2}}$.
\end{theorem} 

In the above theorem, the curvature and torsion were those of the Chern connection. For our purposes, in order to apply later compactness theorems for sequences of manifolds with bounded Hermitian geometry, we need rather bounds for the curvature and torsion of the Bismut connection. But the passage from the Chern connection to any connection on the Gauduchon line can now be obtained from Lemma  \ref{GauduchonEstimates}, applied to the case of fixed $t$. Thus we obtain the following Shi-estimates for the Type IIB flow for any connection $\na^s$ on the Gauduchon line:






\begin{lemma}
	Let $(M,g(t),J)$ be a solution of the Type IIB flow on the interval $[0;\frac{1}{A}]$. 
	As before, we denote by $Rm_s$ and $T_s$ the curvature and the torsion tensors for $\nabla^s$. We denote $Rm_1$ by $Rm$, and $T_1$ by $T$.
	Assume that $|Rm| + |\nabla T| + |T|^2 \leq A$ for any point in $M \times [0;\frac{1}{A}]$. Then, for any $s \in [-1;1]$ there are constants $Q_{s,k}, P_{s,k}$, such that $|(\nabla^s)^k Rm_s| \leq \frac{Q_{s,k}}{t^{k/2}}$, and $|(\nabla^s)^{k+1}T|\leq \frac{P_{s,k}}{t^{k/2}}$.
\end{lemma}





\section{Compactness theorem for the Type IIB flow}
\setcounter{equation}{0}

We give now the proof of Theorem \ref{compactness}.
We break it into several steps.

\medskip
{\it Convergence of the metrics} 

The $C^\infty$ convergence of the Riemannian manifolds $(M_j, g_j)$ to $M_\i, g_\i$, the existence of an exhaustion $U_j$ for $M_\i$ and the existence of diffeomorphisms $\Phi_j: U_j \rightarrow V_j \subset M_j$ all follow from Lemma \ref{CurvatureEstimates} and Hamilton's compactness theorem. 
Indeed, by Lemma \ref{CurvatureEstimates} and Lemma \ref{RiemannianEstiamtes}, we have uniform bounds for the covariant derivatives of Riemannian curvature w.r.t. to the Levi-Civita connection, so we can apply Hamilton's compactness theorem. 

\medskip

{\it Convergence of the complex structures}

Consider the forms $\o_j(X,Y) = g_j(J_jX,Y)$. It is well known that the pair $g_j, \o_j$ defines $J_j$.  Let $\Phi_{j,k}$ will be restriction $\Phi_{j+k}$ on $U_j$. Denote $g_{j,k} := \Phi_{j,k}^*g_{j+k}$. The metrics $g_{j,k}$ converge to $g_{j, \i}$ as $k \rightarrow \i$. It is easy to see that $g_{j, \i}$ is the restriction of $g_\i$ on $U_j$. Due to this, we omit the index $j$ in $g_{j,\i}$ and write $g_\i$ everywhere. We need to show that the forms $\o_{j,k} := \Phi_{j,k}^*\o_{j+k}$ converge uniformly with all derivatives on $U_j$ as $k \rightarrow \i$ to some nondegenerate forms $\o_{j,\infty}$. In order to do this, we need to obtain uniform bounds for $|(\na^\i)^m \o_{j,k}|_{g_\i}$.

Let $v_{j,k} := g_\i - g_{j,k}$ and $A_{j,k} := \na^{g^\i} - \na^{g_{j,k}}$. In a local coordinate chart $(x^1, \dots, x^{2n})$ on $U_j$, the tensor $A_{j,k}$ can be written as follows:
\bea
(A_{j,k})^\a_{\b \g}
=
{1 \over 2}(g^{j,k})^{\a \d}(\na^\i_\b (v_{j,k})_{\d \g} + \na^\i_\g (v_{j,k})_{\b \d} - \na^\i_\d  (v_{j,k})_{\b \g}).
\eea
One can check this formula in the geodesic coordinates for $g_\i$.

Since $g_{j,k}$ converge to $g_\i$ in $C^\i$, the norms $|(\na^{g_\i})^mv_{j,k}|_{g_\i}$ tends to zero for any $m \in \N$. The argument from \cite[Theorem 7.1]{LotayWei15} shows, that for any $m \geq 0$ there is a constant $K_m$, such that $|(\na^{g_\i})^mA_{j,k}|_{g_\i} \leq K_m$.

Now we can obtain bounds for $\o_{j k}$. These bounds follow from the boundedness of the covariant derivatives of $J_j$ (which follows from the Lemma \ref{BoundsOnComplexStructure}), and the bounds for $A_{j,k}$. We have
\bea
\na^{g_\i} \o_{j,k} 
=
 \na^{g_{j,k}}\o_{j,k} + (\na^{g_\i} - \na^{g_{j,k}})\o_{j,k}.  
\eea
Since $g_{j,k}$ converge to $g_\i$ with all derivatives, there is a constant $K_\i$ independent of metrics, such that ${1 \over K_\i}g_{j,k} \leq g_\i \leq K_\i g_{j,k}$ for all $k$ on $U_j$, and
\bea
|\na^{g_\i} \o_{j,k}|_{g_\i} 
\leq 
(K_\i)^{3 \over 2}|\na^{g_j}\o_{j,k}|_{g_{j,k}} + K_0 |\o_{j,k}|_{g_{j,k}}. 
\eea
The quantity $|\na^{g_j}\o_{j,k}|_{g_{j,k}}$ is bounded because the covariant derivatives of $J_j$ are bounded. The quantity $|\o_{j,k}|_{g_{j,k}}$ is a constant.

Next, we obtain bounds for the higher derivatives of $\na_{j,k}$. As in \cite{LotayWei15}, we can write:
\bea
|(\na^{g_\i})^l \o_{j,k}|_{g_\i} 
&\leq&
C\Sigma_{m=0}^l |A_{j,k}|^m_{g_\i} |(\na^{g_{j,k}})^{m-l} \o_{j,k}|_{g_{j,k}} 
\nonumber\\
&&\quad+ 
C\Sigma_{m=1}^{l-1}|(\na^{g_\i})^mA_{j,k}|_{g_\i}|(\na^{g_{j,k}})^{m-l-1}\o_{j,k}|_{g_{j,k}}
\eea
The right hand side is bounded by a constant, because each summand is bounded by a constant. By the Arzela-Ascoli theorem (see \cite[Corollary 9.14]{AndrewsHopper}), we can extract a subsequence, which converge to the form $\o_{j,\i}$ in $C^\i$ on compacts in $U_j$. The form $\o_{j,\i}$ is compatible with $g_\i$ on $U_j$. This follows from the compatibility of $g_{j,k}$ and $\o_{j,k}$. In particular, 
\bea
\o_{j, \i}^n 
= 
\lim_{k\rightarrow\infty}\o_{j, k}^n 
= 
n!\lim_{k\rightarrow\infty}d\mu_{g_{j, k}}  
= 
n! d\mu_{g_{ \i}},
\eea 
hence $\o_{j,\i}$ is nondegenerate.
The pair $(g_{j, \i}, \o_{j,\i})$ defines an almost complex structure $J_{j,\i} := g_{j, \i}^{-1}\o_{j, \i}$, which is the $C^\i$-limit of $J_{j,k} := g_{j, k}^{-1}\o_{j, k}$. In fact, $J_{j, \i}$ is integrable, because all $J_{j,k}$ are integrable, and $\lim_{k\rightarrow\infty}J_{j,k} = J_{j,\i}$ in $C^\i$. 

Let $I_{jk}:U_j \rightarrow U_k$ be the inclusion map for $k \geq j$. We see that $I^*_{jk}\o_{k,l} = \o_{j,l}$, and $I^*_{jk}\o_{k,\i} = \o_{j,\i}$. Since the set of smooth sections of any vector bundle on $M_\i$ is a sheaf, there is a nondegenerate form $\o_\i$ on the whole $M_\i$, such that the restriction of this form on each $U_j$ is $\o_{j,\i}$. Furthermore, the pair $(g_\i, \o_\i)$ defines an almost complex structure $J_\i := g_{ \i}^{-1}\o_{ \i}$ on $M_\i$.

In is not hard to see that the restriction of $J_\i$ on $U_j$ is equal to $J_{j,\i}$, defined by the pairs $(g_{\i}, \o_{j,} )$. Hence, the Nijenhuis tensor of $J_\i$ is zero, and $(M_\i, g_\i, J_\i)$ is a complete Hermitian manifold.

\medskip

{\it Convergence of the holomorphic forms}

By \cite[Lemma 2.5]{Pic}, we know that forms $\psi_j := {\Omega_j \over ||\Omega_j||_{g_j}}$ are parallel w.r.t. Bismut connection $\na^{+,g_j}$. Recall that the Bismut connection on $(n,0)$ forms is written as follows: 
\bea
\na^+_j \psi = \p_j \psi - \G^s_{j s}\psi - T^s_{j s}\psi
\nonumber\\
\na^+_{\bar j} \psi = \p_{\bar j} \psi - g^{s \bar p}\bar{T}_{s \bar j \bar p } \psi
= 
\p_{\bar j} \psi + \bar{\tau}_{\bar j} \psi
\eea
Hence, all derivatives of these forms are bounded, and we again apply the Arzela-Ascoli theorem, and show that the pullbacks $\psi_{j,k}:=\Phi^*_{j+k}\psi_{j+k}$ preconverge to the limit form $\psi_{j,\infty}$ on $U_j$ in the $C^\infty$ sense. Moreover, the form $\psi_{j,\infty}$ will be parallel $(n,0)$-form with respect to the Bismut connection of the pair $(g_\infty, J_\infty)$ and $|\psi_{j,\i}|_{g_\i} = 1$. Indeed, we have the following estimate:
\bea
|\nabla^{+,g_\i}\psi_{j,\infty} - \nabla^{+,g_{j,k}}\psi_{j,k}|_{g_\i} 
\leq
|\nabla^{+, \infty}(\psi_{j,\infty} - \psi_{j,k})|_{g_\i}+|(\nabla^{+,\infty} - \nabla^{+,g_{j,k}})\psi_{j,k}|_{g_\i}. 
\eea
The first summand is bounded by $\varepsilon$, because because metrics, complex structures and $(n,0)$ forms $\psi_{j,k}$ converge to $g_\infty, J_\infty$ and $\psi_{j,\infty}$ in $C^\infty$ sense. The second summand is bounded by $\varepsilon$ by the same reason. As in previous section, there is a globally defined, nondegenerate $(n,0)$-form $\psi_\i$ on the whole $M_\i$, which is parallel w.r.t. the Bismut connection on $(M_\i, g_\i, J_\i)$. 

Let $f_{j,\infty} = \lim_{k\rightarrow\infty}\|\O_{j,k}\|_{g_{j,k}}$. This limit exists and will be the $C^\infty$ function, because the differentials of $\|\O_{j,k}\|$ is the 1-form $\|\O_{j,k}\|\tau(g_{j,k}) = -\|\O_j\|_{g_{j,k}}J_{j,k} d^\dagger_{j,k}\o_{j,k}$ (see \cite[Lemma 4]{PPZ16}; the operator $d^\dagger_{j,k}$ is defined via $\o_{j,k}$). Hence all derivatives of $\|\O_j\|$ are bounded, and we again apply the Arzela-Ascoli theorem.  By our assumption on $\log \|\O_j\|_{g_j}$, this is a $C^\infty$ positive function and $\nabla^{+,g_{j,k}}\|\O_j\|_{g_{j,k}} \rightarrow \na^{+, g_{j,\i}}f_{j,\i} = f_{j,\i}\tau(g_{j,\i})$, and there is a $C^\i$ positive function $f_\i$, defined on the whole $M_\i$, such that its restriction on $U_j$ is $f_{j,\i}$. Consider the form $\O_\i:=f_\infty\psi_\infty$. This form will be a nonvanishing holomorphic form on $(M, J_\infty)$. Indeed, in a local coordinate system on $(M_\i, J_\i)$ we have
\bea
\p_{\bar a} \O_\i
=
\na^+_{\bar a} \O_\i - \tau_{\bar a}\O_\i
=
\na^+_{\bar a} (f_\i \phi_\i) - \tau_{\bar a}f_\i \phi_\i
=
\tau_{\bar a}f_\i \phi_\i - \tau_{\bar a}f_\i \phi_\i
=
0,
\eea
where bar over indices takes w.r.t. $J_\i$, and all $\tau$ is the torsion of the Chern connection on $(M_\i, g_i. J_i)$.
Note that the function $f_\i$ satisfies the following identity: $f_\i = ||\O_\infty||_{g_\i}$.

\medskip

{\it The conformally balanced property of the limit metric}

Finally, the limit metric satisfies the conformally balanced condition. This follows from the convergence of derivatives of metrics and complex structures. Q.E.D.

\bigskip

Clearly, condition (4) in the theorem can be replaced by the condition
$d(\|\Omega_j\|^2_{g_j}\omega^{n-1}_j)=0$, and the conclusion would remain the same, with 
the limiting triple $(g_\infty, J_\infty, \Omega_\infty)$ satisfying the same condition.

\medskip
Finally, we note that a compactness theorem for sequences of K\"ahler manifolds can be found in \cite{_PS06_}.

\subsection{Estimates for sequences of metrics and for $\|\Omega\|$}

In order to prove Theorem \ref{compactnessTypeIIB} we need to estimate the mixed time and spatial derivatives for the metric $g(t)$, which evolves as $\p_t g_{\bar k j} = - \tilde{R}_{\bar k j} - {1 \over 2}T_{\bar k p q}\bar{T}_{j}^{p q}$. We need a fixed Hermitian metric $\hat{g}$, which is independent on time $t$. All geometric quantities  (e.g. the Chern connection, the curvature, etc) related to $\hat{g}$ we will denote with the hat symbol.

Now we are going to prove the following lemma, which is analogous to the Lemma 3.11 from Chow et al.\cite{_CCGGIIKLLN-I_}.

\begin{lemma}\label{BigLemma}
	Let $(M,J,\hat{g})$ be a Hermitian manifold and $K$ be a compact subset of $M$. Assume that $\{g_m(t)\}$ be a sequence of solutions of the flow $\p_t g_{\bar k j} = - \tilde{R}_{\bar k j} - {1 \over 2}T_{\bar k p q}\bar{T}_{j}^{p q}$ on a neighborhood of $K\times[\tau_0;\tau_1]$. Fix some $t_0 \in [\tau_0;\tau_1]$. Assume also that the following conditions holds:

\smallskip
		{\rm (1)} For our $t_0 \in [\tau_0;\tau_1]$, and for any $m$ we have the estimate $$C^{-1}\hat{g} \leq g_m(t_0) \leq C\hat{g}$$
		for the constant $C$ is independent of $k$.
		
		{\rm (2)} There are constants $C^{'}_p>0$, independent of $k$, such that 
		\bea
		|\hat{\nabla}^pg_m(t_0)|_{\hat{g}}
		\leq 
		C^{'}_p
		\eea
		for all $p \geq 1$.

	{\rm (3)} There are constants $C^{''}_p, C^{'''}_p>0$, independent of $k$, such that the covariant derivatives of Chern curvature tensor $Rm$ are bounded 
		\bea
		|(\na^m)^p Rm_m|_{g_m} 
		\leq 
		C^{''}_p
		\eea
		Moreover, all covariant derivatives of the torsion $T_m$ are bounded as well
		\bea
		|\na^pT_m|_{g_m} 
		\leq 
		C^{'''}_p.
		\eea
		Here $\na^m$ denotes the Chern connection of $g_m(t)$.

\smallskip
	
	Then the following holds
	
	\smallskip
	
		{\rm (1)} There are constant $N$ independent of $m$, such that  \bea
		N^{-1}\hat{g} 
		\leq 
		g_m(t) 
		\leq 
		N\hat{g}
		\eea
		
	{\rm (2)} There are constants $N^{'}_p>0$, independent of $k$, such that 
		\bea
		|\hat{\nabla}^pg_m(t)|_{\hat{g}}
		\leq 
		N^{'}_p
		\eea
		for all $p \geq 1$.
		
		{\rm (3)} There are constants $N^{''}_{p,q}>0$, such that 
		\bea
		|\frac{\partial^q}{\partial t^q}\hat{\nabla}^pg_m(t)|_{\hat{g}}
		\leq 
		N^{''}_{p,q}
		\eea

\end{lemma}

\noindent{\em Proof.}
{\rm (1)}  For any vector $V$ and for any solution $g$ of the Type IIB flow we have the following identity
	\bea
	\partial_t g(V,V) 
	=
	- \tilde{R}_{j\bar{k}}V^j\overline{V}^k - {1 \over 2}(T \circ \bar T)_{j\bar{k}}V^j\overline{V}^k.
	\eea
	
	From the definition we see that $|\tilde{R}_{j\bar{k}}V^j\overline{V}^k| \leq \sqrt{n}|Rm|g(V,V)$, and the norm of $T \circ \bar T$ is bounded.
	
	Hence, by the standard trick, we can write the following:
	\bea
	C|t_1 - t_0| \geq \int_{t_0}^{t_1}|\partial_t \log(g(V,V))|dt 
	\geq 
	|\log\frac{g(t_1)(V,V)}{g(t_0)(V,V)}|
	\eea
	for some constant $C$, which depends on the dimension, and the bounds for $T$ and $Rm$. Now we can take the exponent and obtain the desired bound with $N=e^{C(t_1 - t_0)}$
	
	\smallskip
	{(2)} Throughout this section $\na^m$ and $\G_m$ denote the Chern connection with respect to $g_m$.
	
	We want to obtain bounds for $|\hat{\na}^pg_m(t)|_{\hat{g}}$ for any $p$. First of all (see \cite{FPPZa} for details), $h^{-1}_m\hat{\na}h_m = \na^m h_m h^{-1}_m = \na^m - \hat{\na} = g^{-1}\hat{\na}g$, where $h_m$ is a relative endomorphism $(h_m)^a_b = \hat{g}^{a\bar{s}}(g_m)_{b\bar{s}}$. Next, since we proved the equivalence for $\hat{g}$ and $g_m$, the norms on tensors, induced by these metrics are equivalent (see [Lemma 3.13]\cite{_CCGGIIKLLN-I_} for the proof). 
	
	We have the following identity:
	\bea
	\p_t \hat{\na}g_m
	=
	\hat{\na} \dot{g}
	= 
	-\hat{\na} (\tilde{Ric}_m + (T \circ \bar T)_m)
	=
	-\na(\tilde{Ric}_m + (T \circ \bar T)_m) +
	\nonumber\\ 
	(\na - \hat{\na})(\tilde{Ric}_m + (T \circ \bar T)_m)
	\eea
	
	The first summand is bounded by our assumptions. The second summand is bounded by constant times $|\hat{\na}^pg_m(t)|_{\hat{g}}$. Altogether, we obtain the inequality
	\bea
	|\hat{\na}^pg_m(t)|_{\hat{g}}
	\leq
	C(1 + |\hat{\na}^pg_m(t)|_{\hat{g}}),
	\eea
	where $C$ is the maximum of bounds for $\tilde{Ric}_m + (T \circ \bar T)_m$ and for $\na^m(\tilde{Ric}_m + (T \circ \bar T)_m)$.
	
	For $p>1$ we proceed by induction.  First of all, for any $p \in \N$, by the same way as in \cite[Lemma 8.6]{AndrewsHopper}, one can show that
	
	\bea
	\p_t \hat{\nabla}^p g_m
	=
	\sum_{k_1 + \dots + k_q = p} (\na^m)^{k_1}(Rm_m + (T * \bar T)_m) * \hat{\na}^{k_2}g_m * \dots * \hat{\na}^{k_q}g_m.
	\eea
	
	Suppose we have the estimate for $p-1$. Then, taking the time derivative of $|\na^p g_m|_{\hat{g}}$ and incorporate bounds for the curvature, the torsion and for the lower derivatives, we obtain:
	\bea
	\p_t |\na^p g_m|_{\hat{g}}
	\leq
	C(1 + |\na^p g_m|_{\hat{g}}),
	\eea
	where $C>0$ depends only on bounds for the curvature, the torsion and for the lower derivatives. Hence,
	\bea
	\log(|\na^p g_m|_{\hat{g}}+1)
	\leq
	Ct 
	\eea
	for some constant $D$. Now we take exponent and obtain the desired bound.
	
	\smallskip
	{\rm (3)} We have 
	\bea
	|\p^q_t\hat{\na}^pg_m(t)|_{\hat{g}} 
	=
	|\hat{\na}^p\p^q_tg_m(t)|_{\hat{g}}
	=
	|\hat{\na}^p\p^{q-1}_t\dot{g}_m(t)|_{\hat{g}}.
	\eea
	
	The Laplacian formulation of the Type IIB flow (see \ref{laplacian}) allows us write time derivatives of the metric $g_m$ via the spatial derivatives. Iterating this process, we are able to estimate all mixed derivatives via spatial derivatives. However, we proved before the estimates for spatial derivatives.  Q.E.D.
	
	\smallskip

We can derive a useful consequence of this lemma.

\begin{lemma}\label{BigCorollary}
	Let $(M,\hat{g},J, \O)$ be a Hermitian Calabi-Yau manifold and $K \subset M$ is compact. Assume there are complex structures $\{J_j\}$ and holomorphic forms $\{\O_j\} $ on $M$, which converge in $C^\i$ to $J$ and $\O$ respectfully. Suppose $g_j(t)$ are the solutions of the Type IIB flow w.r.t. $J_j$ and $\O_j$ on $K \times [\tau_0; \tau_1]$. Suppose all $g_j(t)$ satisfy the assumptions of Lemma \ref{BigLemma}. Then the results of Lemma \ref{BigLemma} also holds for $g_j(t)$, i.e. there are constants $N,N_{l,m}$, such that $|\p^l_t \hat{\na}^m g_j(t)|_{\hat{g}} \leq N_{l,m}$ and $N^{-1}\hat{g} \leq g_j \leq N\hat{g}$.
\end{lemma}

The difference between this Lemma and Lemma \ref{BigLemma} is that here, the metrics $g_j(t)$ are Hermitian with respect to the {\em different} complex structures on $K$.

\medskip

\noindent{\em Proof.}
Define the sequence of the metrics $\hat{g}_j$ by the formula $\hat{g}_j := {1 \over 2}(\hat{g}(X,Y) + \hat{g}(J_jX,J_jY))$. The metrics $\hat{g}_j$ converge to $\hat{g}$ in $C^\i$ because $J_j$ converge to $J$ in $C^\i$. Without loss of generality, we can assume that for any $j$ there is a constant $P$ independent on $j$, such that $P^{-1}\hat{g} \leq \hat{g}_j \leq P \hat{g}$. Applying the estimates from Part (1) of the previous lemma, we obtain the uniform equivalence between $g_j(t)$ and $\hat{g}$.

Denote by $\hat{\G}_j$ and $\hat{\G}$ the Chern connections of $(\hat{g}_j,J_j)$ and $(\hat{g},J)$ respectfully. The difference $|\hat{\G}_j - \hat{\G}|_{\hat{g}}$ between Chern connection $\hat{\G}_j$ of $(\hat{g}_j, J_j)$ and $\hat{\G}$ of the pair $(\hat{g},J)$ tends to zero as $j \rightarrow \i$. Let $\hat{S}_j(X,Y,Z) := d\hat{\o}_j(J_jX,Y,Z)$ and $\hat{S}$ has the same meaning for $(\hat{g}, J)$. The difference of the Chern connections can be written as follows  
\bea
\hat{\G}_j - \hat{\G}
=
(\na^{\hat{g}_j} - \na^{\hat{g}}) + (\hat{g}^{-1}_jS_j - \hat{g}^{-1}\hat{S}).
\eea
The quantity $|\na^{\hat{g}_j} - \na^{\hat{g}}|_{\hat{g}}$ tends to zero as $j$ tends to infinity due to $C^\i$ convergence of $\hat{g}_j$ and $J_j$. The norm of the second summand is again tends to zero by the same reason.

The same is true for $|\hat{\na}^l(\hat{\Gamma}_j - \hat{\Gamma})|_{\hat{g}}$ for any $l \in \N$. Hence $|\hat{\na}^l(\hat{\Gamma}_j - \hat{\Gamma})|_{\hat{g}}$ is bounded by a constant, independent on $j$. Now we apply the previous lemma to $g_j(t)$ and $\hat{g}_j$, and use the fact that $|\hat{\na}^l(\Gamma_j - \hat{\Gamma})|_{\hat{g}} \leq |\hat{\na}^l(\Gamma_j - \hat{\Gamma}_j)|_{\hat{g}} + |\hat{\na}^l(\hat{\Gamma}_j - \hat{\Gamma})|_{\hat{g}}$. Hence, we obtain bounds for $|\hat{\na}^lg_j(t)|_{\hat{g}}$ and $|\p^l \hat{\na}^m g_j(t)|_{\hat{g}}$ with respect to $\hat{g}$. Q.E.D.

\medskip

Finally, we need to study the evolution of $\log\|\O\|^2$ along the flow.

\begin{lemma}\label{DilatonBounds}
	Let $(M, g(t),J, \O)$ be a compact conformally balanced Calabi-Yau manifold, and $g(t)$ evolves according to the Type IIB flow: $\p_t g_{\bar k j} = - \tilde{R}_{\bar k j} - {1 \over 2}T_{\bar k p q}\bar{T}_{j}^{p q}$. Assume there is a constant $C>0$, such that $|Rm|_g + |\na T|_g + |T|^2_g \leq C$. Then there are constants $C_1, C_2$ and $C_3$ such that the following holds:
	\bea
	C_1 \leq \log\|\O\|^2_g(t) \leq C_2t + C_3. 
	\eea
\end{lemma}

\noindent{\em Proof.}
By a straightforward computation, we can show that $\tilde{R} = g^{j \bar k} \tilde{R}_{\bar k j} = g^{j \bar k} \p_{\bar k} \p_j \log\|\O\|^2 = \Delta \log\|\O\|^2$ and, 
\bea
(\p_t - \Delta)\log\|\O\|^2 
= 
{1 \over 2}|T|^2_g
\geq 0
\eea
(see \cite[Section 3.1]{FPPZa}). By the maximum principle, $\log\|\O\|^2$ achieves minimum at $t=0$, hence it bounded below.
On the other hand, by the assumption on $Rm$ and $T$, there is constant $C_2$, such that the following inequality holds:
\bea
\p_t \log\|\O|\|^2 = \tilde{R} + {1 \over 2}|T|^2_g \leq C_2.
\eea
Now we can integrate it along $t$, and obtain the bound from above. Q.E.D.

\subsection{Proof of Theorem \ref{compactnessTypeIIB}}

We proceed as in the proof of Hamilton's compactness theorem (see \cite{Hamilton1, LotayWei15}).

By Theorem \ref{compactness}, Shi-type estimates from \cite{PPZ16}, \cite{_FP19_}, and by Lemma \ref{DilatonBounds}, we can guarantee that the family $\{(M_j, g_j(\tau_0), J_j, \O_j, x_j)\}_{j=1}^{\infty}$ converge to the manifold $(M_\infty, \hat{g}, J_\infty, \O_\i,  x_\i)$. Note that Shi-type estimates, together with Lemma \ref{BoundsOnComplexStructure} imply the bounds on the covariant derivatives of complex structures. Moreover, since the bounds on torsion does not depend on $g_j$, the bounds on the covariant derivatives of $J_j$ do not depend on $g_j$ either. 

\smallskip

Now we apply Lemma \ref{BigLemma} and Lemma \ref{BigCorollary} for $\hat{g}$ and the sequence of metrics $\tilde{g}_j:=\Phi^*_jg_j(t)$ on a compact set $K \subset U_j \subset M_\i$, and obtain uniform $C^\i$ control on these metrics. The Shi-type estimates from \cite{PPZ16} ensure that the conditions of Lemma \ref{BigLemma} are satisfied. Hence, we can consider the sequence $\tilde{g}_j(t) + dt^2$, which will converge on $K \times [\tau_0 ; \tau_1]$ to some metric $g_\infty(t) + dt^2$ (by the Arzela-Ascoli theorem). The limit $g_\infty(t)$ will be the solution to the Type IIB flow with the initial metric $\hat{g}$, i.e. $g_\infty(\tau_0) = \hat{g}$ with respect to the complex structure $J_\i$ and the holomorphic volume form $\O_\i$. This follows from the $C^\i$ convergence of metrics to the limit metric. By the diagonalization argument, there is a subsequence, which converges on the whole $M_\i \times [\tau_0; \tau_1]$  
Q.E.D.

\hfill

\section{Singularity models for the Type IIB flow}
\setcounter{equation}{0}

Let $\o(t)$ be a solution for the Type IIB flow on a complex manifold on some time interval $[0,T)$. For each $t\in [0,T)$, define the function $f(t)$ by 
\bea
\label{f}
f(t):=\sup_M(|Rm|^2 + |\nabla T|^2 + |T|^4).
\eea

\begin{definition}[Singularity types]
	Let $(M, g(t), J, \Omega)$ be a solution of the Type IIB flow on a compact Hermitian Calabi-Yau manifold $M$. Let $f(t)=\sup_M(|Rm|^2 + |\nabla T|^2 + |T|^4)$ as above, and $T$ is the maximal time for which flow exists. We define the following type of singularities
	
		{\rm [Type I]} $\sup_{[0;T]}(T-t)f(t) = C < \infty$, and $T<\infty$	
		
		{\rm [Type IIa]} $\sup_{[0;T]}(T-t)f(t) = \infty$, and $T<\infty$
		
		{\rm [Type IIb] }$\sup_{[0;T]}tf(t) = \infty$, $T=\infty$
		
		{\rm [Type III]} $\sup_{[0;T]}tf(t) = C < \infty$, $T = \infty$
	
\end{definition}

We observe that we could 
also try to define singularities in terms of $\|\O\|^2$, as suggested by examples (see \cite[Section 2.3.2]{Pic} for instance). However, as Lemma \ref{DilatonBounds} shows, the bounds on $f(t)$ imply bounds on $\|\Omega\|$. 
\hfill
\begin{definition}
\label{injestimate}
	We say that $(M, g(t), J, \Omega)$ satisfies the injectivity radius estimate if the injectivity radius satisfy the following inequality
	\bea
	\mathrm{inj}(M, g(t)) \geq {c \over \sqrt{f(t)}}
	\eea
	for some positive constant $c$.
\end{definition}

We give now the proof of Theorem \ref{models}.
Denote $F(x,t) := \sqrt{|Rm|^2 + |\nabla T|^2 + |T|^4}$. Recall that $f(t) = \sup_M\sqrt{|Rm|^2 + |\nabla T|^2 + |T|^4} = \sup_{x \in M} F(x,t)$. Throughout the proof of the theorem we will use the sequence of rescaled Hermitian metrics $g_j(t)$, defined as follows:
\bea
g_j(t) := C_jg(t_j + C_j^{-1}t),
\eea
where the constants $C_j$ will be defined later. We also rescale the holomorphic volume form
\bea
\Omega_j :=  C_j^{-n \over 2}\Omega.
\eea

\medskip

\noindent{\em The case of Type I singularity}

Let $C := \sup_t f(t)$. . Pick a sequence $(x_j, t_j) \in M \times [0;T]$, such that $F(x_j,t_j) \rightarrow C$. Let $C_j := F(x_j,t_j)$, and $(T-t_j)F(x_j,t_j) \geq c >0$.
Consider the rescaled Hermitian metrics $g_j(t)$ and forms $\Omega_j$. These metrics are again the solutions of the Type IIB flow on the interval $(-C_jt_j; C_j(T-t_j))$. The curvature and torsion tensors of the metrics $g_j$ satisfy the following bound:
\bea
\sqrt{|Rm|_{g_j}^2 + |\nabla T|_{g_j}^2 + |T|_{g_j}^4} \leq C_j^{-1}{C \over T - (t_j + C_j^{-1}t)} = {C \over C_j(T - t_j) - t } \leq {C \over c - t}.  
\eea
By the compactness theorem for Type IIB flow, the sequence $(M, g_j(t), J, \Omega_j, x_j)$ converge to the limit manifold $(M_\infty, g_\infty(t), J, \Omega, \infty)$ with $|Rm|_{g_\infty}^2 + |\nabla T|_{g_\infty}^2 + |T|_{g_\infty}^4 \leq {C \over c -t}$. The case of Type I singularity is proved.

\medskip

\noindent{\em The case of Type IIa singularity}

Pick a monotone increasing sequence ${T_j}$, such that $T_j \rightarrow T$ and let $C_j = \sup_{t<T_j}f(t)$. Next, we pick a sequence of points $(x_j,t_j)$, such that $C^2_j = (T-t_j)F(x_j,t_j)$. Now we take a sequence of rescaled metrics $g_j(t) := C_jg(t_j + C_j^{-1}t)$ and the sequence of pointed manifolds $(M, g_j(t), J, \Omega_j, x_j)$ again converge to the limit manifold $(M, g_\infty(t), J, \Omega, x_\infty)$. The limit metric $g_\infty(t)$ is the solution of the Type IIB flow on the whole line $(-\infty,\infty)$, as $C_j(T-t_j) \rightarrow \infty$. The bound for the curvature is the following:
\bea
\sqrt{|Rm|_{g_j}^2 + |\nabla T|_{g_j}^2 + |T|_{g_j}^4} \leq  {C_j \over C_j(T - t_j)} \rightarrow 1.
\eea
Hence, $f(t) \leq 1$.

\medskip

\noindent{\em The case of Type IIb singularity}
Pick a monotone sequence $T_j \rightarrow \infty$ and a sequence of points $(x_j, t_j), t_j \leq T_j$, such that:
$t_j(T-t_j)F(x_j,t_j) = \sup t(T_j - t)F(x,t)$. Let $C_j$ be as above, and $A_j = t_jC_j$

\bea
\sqrt{|Rm|_{g_j}^2 + |\nabla T|_{g_j}^2 + |T|_{g_j}^4} \leq C_j^{-1} {A_j \over t_j + C_j^{-1}t} 
=
{A_j \over A_j + t} 
\rightarrow 
1
\eea

We again apply the compactness theorem and have a limit with the desirable property.

\hfill

\noindent{\em The case of Type III singularity}

Let $A = \lim \sup tf(t) < + \infty$. This number is positive. Indeed, if $tf(t)< \varepsilon$, then the diameter $L$ satisfies the following inequality:
\bea
{\partial L\over \partial t} \leq {C\varepsilon L \over t},
\eea
because the $f(t)$ is bounded, hence $\tilde{R}_{\bar k j} + {1 \over 2}T_{\bar k p q}\bar{T}_j^{p q}$ is bounded. Hence, $\log L \leq C\varepsilon \log t$. Hence, $L \leq t^{C\varepsilon}$. Take $\varepsilon$ smaller than ${1 \over 2C}$. This choice implies $L^2f(t) \rightarrow 0$. On the other hand, $L \geq \mathrm{inj}(M, g(t)) \geq {C \over \sqrt{f(t)}}$ by the assumption of the theorem. This inequality is equivalent to the inequality $L^2f(t) \geq C$. This is a contradiction, hence $A>0$.

\smallskip

Let $C= \sup f(t)$, and $C_j, (x_j,t_j)$ as above. The metrics $g_j(t)$ are solutions to the Type IIB flow on the interval $[-C_jt_j; +\infty)$ and, by our assumptions $C_jt_j \rightarrow A$.  We have the following inequality:
\bea
\sqrt{|Rm|_{g_j}^2 + |\nabla T|_{g_j}^2 + |T|_{g_j}^4} 
\leq 
C^{-1}_j{A \over t_j +C_j^{-1}t} 
=
{A \over C_jt_j + t } 
\rightarrow 
{A \over A + t} .
\eea	
Q.E.D.

\bigskip

\noindent Department of Mathematics, Columbia University, New York, NY 10027 USA

\noindent also

\noindent Laboratory of Algebraic Geometry,
Faculty of Mathematics, HSE University, 6 Usacheva Str. Moscow, Russia

\noindent nklemyatin@math.columbia.edu

\end{document}